\documentclass[12pt]{article}
\usepackage{theorem,amsfonts}

\newtheorem{theorem}{Theorem}[section]
\newtheorem{proposition}{Proposition}[section]
\newtheorem{lemma}{Lemma}[section]
\newtheorem{corollary}{Corollary}[section]
\theorembodyfont{\upshape}
\newtheorem{definition}{Definition}[section]
\newtheorem{remark}{Remark}[section]
\newtheorem{proof}{Proof}

\newtheorem{acknowledgement}{Acknowledgement}

\newcommand{\bt}{\begin{theorem}}
\newcommand{\et}{\end{theorem}}
\newcommand{\bl}{\begin{lemma}}
\newcommand{\el}{\end{lemma}}
\newcommand{\bp}{\begin{proposition}}
\newcommand{\ep}{\end{proposition}}
\newcommand{\bo}{\begin{proof}}
\newcommand{\eo}{\end{proof}}
\newcommand{\br}{\begin{remark}}
\newcommand{\er}{\end{remark}}
\newcommand{\bc}{\begin{corollary}}
\newcommand{\ec}{\end{corollary}}
\newcommand{\bd}{\begin{definition}}
\newcommand{\ed}{\end{definition}}
\newcommand{\be}{\begin{enumerate}}
\newcommand{\ee}{\end{enumerate}}

\title{Absolute continuity of autophage measures on finite-dimensional vector
spaces}
\author{C. R. E. Raja } 

\date{}

\begin{document}
\maketitle

\let\epsi=\epsilon
\let\vepsi=\varepsilon
\let\lam=\lambda
\let\Lam=\Lambda 
\let\ap=\alpha
\let\vp=\varphi
\let\ra=\rightarrow
\let\Ra=\Rightarrow 
\let\LRa=\Leftrightarrow
\let\Llra=\Longleftrightarrow
\let\Lla=\Longleftarrow
\let\lra=\longrightarrow
\let\Lra=\Longrightarrow
\let\ba=\beta
\let\ga=\gamma
\let\Ga=\Gamma
\let\un=\upsilon

\begin{abstract}
We consider a class of measures called autophage which was introduced and studied by
Szekely for measures on the real line.  We show that the autophage measures on 
finite-dimensional vector spaces over real or ${\mathbb Q}_p$ are infinitely 
divisible without idempotent factors and are absolutely continuous with bounded 
continuous density.  We also show that certain semistable measures on such vector 
spaces are absolutely continuous.  
\end{abstract}

\begin{section}{Introduction and Preliminaries}

Let $G$ be a locally compact second countable group.  Let ${\cal P}(G)$ be the space
of Borel probability measures on $G$, equipped with the weak topology.  For any
$\mu$ and $\lam \in {\cal P}(G)$, $\mu *\lam$ denotes the convolution of $\lam$
and $\mu$.  For $\mu \in {\cal P}(G)$, $\check \mu$ is defined by $\check \mu (B) =
\mu (B^{-1})$ for any Borel subset $B$ of $G$ and for any integer $n \geq 1$, 
$\mu ^n$ denotes the $n$th-convolution of power of $\mu$.  We say that a measure
$\mu \in {\cal P}(G)$ is {\it continuous} if $\mu (\{x \}) =0$ for all $x \in G$ and
{\it absolutely continuous} if $\mu$ is absolutely continuous with respect to a Haar
measure on $G$.   A {\it continuous convolution semigroup} $(\mu _t)$ in 
${\cal P}(G)$ is a continuous homomorphism $t\mapsto \mu _t$ from the additive
semigroup ${\mathbb R}^+ \cup {0}$ into ${\cal P}(G)$.  

Let Aut$~(G)$ denotes the group of all continuous automorphisms of $G$.  Then for
any $T \in {\rm Aut}~(G)$ and $\lam \in {\cal P}(G)$, $T(\lam)$ is defined by 
$T(\lam ) (B) = \lam (T^{-1}(B))$ for any Borel subset $B$ of $G$.  An automorphism  
$T\in {\rm Aut}~(G)$ is called {\it contraction} if $T^nx \ra e$ for all 
$x \in  G$.  In case, $G$ is a finite-dimensional vector space over real or 
${\mathbb Q}_p$, we say that $T \in {\rm Aut}~(G)$ is a {\it strict-contraction} if
$||T||<1$, where $|| \cdot ||$ is the operator norm. 

We say that a $\lam \in {\cal P}(G)$ is a factor of a $\mu \in {\cal P}(G)$ if there
exists a $\nu \in {\cal P}(G)$ such that $\mu = T(\lam )*\nu$ and $\nu$ is called
co-factor.  Here we wish to study the measures $\mu$ for which there exists
commuting contractions $T$ and $S$ such that $T(\mu )$ is a factor of $\mu$ with
$S(\mu )$ as the corresponding co-factor.  This type of measures on the real line 
was introduced and studied by Szekely and Szekely calls such measures autophage
(see \cite {Sz}).  He also proved that nondegenerate autophage measures are
infinitely divisible and absolutely continuous with a bounded continuous density.  
The main part of his proof was to prove that the Fourier transform of such measures
are integrable.  This implies by Fourier inversion formula that such measures are
absolutely continuous with continuous bounded density.  Here we prove that autophage
measures on finite dimensional vector spaces over real or ${\mathbb Q}_p$ are
infinitely divisible and absolutely continuous.  For a given infinitely
divisible probability measure it is of considerable interest to know when it is
absolutely continuous.  In this point of view these results are of considerable
interest.  We also show that certain semistable
measures on vector spaces over real and ${\mathbb Q}_p$ are absolutely continuous
with bounded continuous density: see \cite {Lu} for results on the absolute
continuity of semistable measures on  ${\mathbb R}^N$.  

\end{section}

\begin{section}{Menu-card of measures}
\bd
We say that a measure $\nu \in {\cal P}(G)$ belongs to the {\it menu-card} of a
measure $\mu \in {\cal P}(G)$ if $\mu = T(\mu) *\nu$ for some contraction $T$ of
$G$.  
\ed

\br
It is easy to see that a menu-card of a measure $\mu$ is non-empty if and only if
$\mu$ is strongly $T$-decomposable for some contraction $T$: a measure $\mu \in 
{\cal P}(G)$ is strongly $T$-decomposable for a $T\in {\rm Aut}~(G)$ 
if $T(\mu )$ is a factor of $\mu$ and $T^n (\mu ) \ra \delta _e$ 
(see \cite {Si} and \cite {RS} for results regarding strongly $T$-decomposable
measures).
\er

We now prove that the nondegenerate strongly $T$-decomposable measures on real and 
$p$-adic unipotent groups, studied in \cite {RS} are continuous which implies the
continuity of measures with nontrivial menu-card.  By an algebraic group over real
or ${\mathbb Q}_p$ we mean the group of real or ${\mathbb Q}_p$ points of an
algebraic group defined over real or ${\mathbb Q}_p$ respectively.  

\bp\label{prp0}
Let $G$ be a unipotent algebraic group over real or ${\mathbb Q}_p$.  Then 
nondegenerate strongly $T$-decomposable measures on $G$ are diffuse, that is,
continuous.  In particular, if the menu-card of a measure $\mu \in {\cal P}(G)$
contains a nondegenerate measure then $\mu$ is continuous.
\ep

\bo
Let $T \in {\rm Aut}~(G)$ be such that $\mu = T(\mu ) *\nu$ for some $\nu \in {\cal
P}(G)$ and $T^n (\mu ) \ra \delta _e$.  Suppose $\mu$ is not continuous.  We now claim 
that $\mu$ degenerates.   Let $Z$ be the center of $G$.  Let $G_0 = G$
and for $i \geq1$, $G_{i}= [G, G_{i-1}]$.  There exists a $k\geq 0$ such that 
$G_k \not = (e)$ but $G_{k+1} =(e) $.  Proof is based on induction on $k$.  
Suppose $k=0$, then $G$ is abelian.  Since $G$ has no elements of finite order, 
by 3.5 of \cite {Si}, $\nu$ degenerates and hence $\mu$ degenerates.  
Suppose $k >0$.  Let $G' = G/Z$ and $\mu '$ and $\nu'$ be the projection of $\mu$
and $\nu$ under the canonical map from $G$ to $G'$.  Let $T'$ be the factor 
automorphism of $G'$ corresponding to $T$.  We then have $\mu '= T'\mu ' *\nu'$.  
It may be easily seen that $\mu '$ is also not continuous and hence by induction
$\mu '$ as well as $\nu '$ degenerates.  
Now, let $\lam = \mu *\check \mu $ and $\rho = \nu *\check \nu$.  Then 
$\lam = T(\mu ) * \nu *\check \nu *T(\check \mu ) = T(\lam )*\rho$ since $\nu*\check
\nu$ is supported on the center of $G$.  Also, $\lam$ and $\rho$ are supported on
the center of $G$ and the center of $G$ is $T$-invariant.  Since $\mu$ is not
continuous, there exists a $x_0 \in G$ such that $\mu (\{ x_0 \}) >0$.  Now, 
$\lam (\{ e \}) = \mu *\check \mu (\{ e \}) \geq \mu (\{ x_0 \})^2 >0$.  
Thus, $\lam $ is also not continuous.  Thus, we infer from the abelian case that 
$\lam$ degenerates and so does $\mu$.  
\eo

\end{section}

\begin{section}{Autophage measure}

\bd\label{defn}
A measure $\mu \in {\cal P}(G)$ is called {\it autophage} if there exists
contractions $T$ and $S$ of $G$ such that $\mu = T(\mu ) *S(\mu )$ and $TS=ST$.  
\ed

In \cite {Sz} and \cite {RuSz} autophage measures were introduced on the real line
$\mathbb R$.  We recall the definition of autophage on $\mathbb R$, $\mu \in {\cal
P}(\mathbb R)$ is autophage if $F(x) = F(ax) F(bx)$ for some $a, b>1$ where $F$ is
the distribution of $\mu$.  Our generalization in Definition \ref{defn} for measures
on general groups, we replace the homothesies by two commuting contractions.  If we
had chosen two elements of a one-parameter semigroup of contractions in place of
homothesies, our notion would be very restrictive.  Because in the $p$-adic 
situation, there are no one-parameter groups but there are many autophage measures 
in the sense of Definition \ref{defn}.  Thus, our choice of commuting contractions 
in Definition \ref{defn} yields more autophage measures.

It may be easily seen that stable measures on finite-dimensional vector spaces
are autophage.  More generally, measures that are stable with respect to a
commutative group of contractive automorphisms are autophage.  We recall that a
measure $\mu$ in ${\cal P}(G)$ is said to be $(T, n)$-semistable for some 
$n \in \mathbb N$ and for an automorphism $T$ if $T (\mu ^n) = \mu $.   
It may be easily seen that $(T,2)$-semistable measure is autophage when $T$ is 
a contractive automorphism.  

We now directly show that Gaussian measures on ${\mathbb R}^n$ are
autophage.  Let $\mu$ be a symmetric Gaussian measure on ${\mathbb R}^n$.
Then $\hat \mu (x) = \exp -<Px,x>$ where $P$ 
is a positive definite matrix.  Let $T$ be any strict-contraction on ${\mathbb R}^n$ 
commuting with $P$ such that $T$ is invertible and symmetric (we may choose $T$ to 
be a non-zero homothesy).  Then for any $x \in {\mathbb R}^n$, 
$0 <||x|| ^2 - ||T(x)||^2 = <I-T^2 (x), x>$, that is,  $I-T^2$ is positive
definite.  Let $S$ be the unique square root of $I-T^2$.  Then $S$ is invertible 
and commutes with $T$.  
Since $T$ commutes with $P$, $S$ also commutes with $P$.  This implies that 
$TPT^* + SPS^* = T^2 P +S^2 P = (T^2 +S^2)P =P$.  Thus, $\mu = T(\mu ) *S(\mu)$ 
where $T$ and $S$ are mutually commuting strict-contractions.

\end{section}

\begin{section}{Infinite divisibility}

For any two contractions $T$ and $S$ on $G$, let $sg (T,S)$ be the semigroup
generated by $T$ and $S$.  For any $\ap \in sg (T,S)$, define $l(\ap )$ to be the
sum of number of $T$ and $S$ appear in the expression of $\ap$.  Then we have 

\bl\label{lem1}
Let $V$ be a finite-dimensional vector space over real or ${\mathbb Q}_p$.  Let $T$ 
and $S$ be strict-contractions on $V$.  Then  $\ap (x) \ra e$ for all $x \in V$ as 
$l(\ap ) \ra \infty$.
\el

\bo
Since $T$ and $S$ are strict-contractions, for a given $\epsilon >0$, there exists a
$k \geq 1$, $||T||^n<1$ and $||S||^n<1$ for all $n \geq k$.  For $\ap \in sg(T,S)$
with $l(\ap )=n$, we have $||\ap || \leq ||T||^{n-i} ||S||^i$ for some $0\leq i 
\leq n$.  If $n=l(\ap ) \geq 2k$, then $||\ap || < \epsi$.  Thus, $||\ap || \ra 0$
as $l(\ap ) \ra \infty$.
\eo

\bp\label{pr1}
Let $V$ be a finite-dimensional vector space over real or ${\mathbb Q}_p$.  Let
$T_1,T_2, \cdots , T_m$ be strict-contractions on $V$.  Suppose $\mu$ is a 
probability measure such that $\mu = T_1(\mu ) * T_2 (\mu )*\cdots *T_m(\mu)$. Then 
$\mu$ is infinitely divisible without idempotent factors.  
\ep

\bo
We treat the case $m=2$ and the general case may be done in a similar fashion.  
Let $T$ and $S$ be strict-contractions on $G$ such that $\mu = T(\mu )*S(\mu )$.  We
first claim that $\mu$ is infinitely divisible.  Iteratively we get that 
$$\begin{array}{rcl}
\mu & = & T(\mu )S(\mu ) \\
&= & T^2 (\mu) *TS(\mu ) *ST(\mu ) * S^2 (\mu ) \\
&= & T^3(\mu ) * T^2S(\mu ) * TST(\mu ) * TS^2 (\mu ) * ST^2 (\mu ) *STS
(\mu) *S^2T(\mu)*S^3 (\mu ) \\
&= & ........ \\ 
&= & \prod _\ap \ap (\mu) \\
\end{array}$$  
where the product varies over all $\ap$ in $\{ \ap \in sg(T, S) \mid l(\ap )=n \}$.  
Thus, we have a triangular system of probability measure where sum of the row
product is $\mu$ and the nth row is 
$$T^n (\mu ), T^{n-1}S(\mu ), \cdots , S^{n-1} T(\mu ), S^n (\mu ). $$  
Since the group
is abelian, it is enough to claim that the triangular system is infinitesimal.  
That is, we should prove that $\ap (\mu ) \ra e$ as $l(\ap ) \ra \infty$ which
follows from Lemma \ref{lem1} and dominated convergence theorem.

We now claim that $\mu$ has no idempotent factors.  We may assume that $V$ is a
$p$-adic vector space.  Let $\omega _K$ be the maximal idempotent factor of
$\mu$.  By a Theorem of Parthasarathy and Sazanov (see Theorem 7.2 of \cite {Par}), 
$\mu = \omega _K * \mu_1$ where $\mu_1$ is an
infinitely divisible measure without idempotent factors.  Since $T$ and $S$ are
automorphisms, $T(K)$ and $S(K)$ are the maximal idempotent factors of $T(\mu )$ and
$S(\mu )$ respectively and $T(\mu _1)$ and $S(\mu _1)$ are infinitely divisible
measures without idempotent factors.  We also get that 
$\omega _K * \mu _1 = \omega _{T(K)S(K)}* T(\mu _1)* S(\mu _1)$.  Since $T(\mu _1)$
and $S(\mu _1)$ are infinitely divisible without idempotent factors, by Theorem 4.2
of \cite {PRV}, the Fourier
transform of $T(\mu _1)$ and $S(\mu _1)$ do not vanish .  Let $\chi$ be any
character of $V$.  If $\hat \omega _K (\chi ) =0$, then $\hat \omega
_{T(K)S(K)}(\chi ) =0$.  If $\hat \omega _{T(K)S(K)} (\chi )=0$, then since $\hat
\mu _1 (\chi )\not =0$, we have $\hat \omega _K(\chi ) =0$.  Thus, $\omega _K =
\omega _{T(K)S(K)}$, that is, $K= T(K)S(K)$.  Iteratively we prove that 
$\prod _{l(\ap ) =n} \ap (K) = K$ for any $n \geq 1$.  Since $\ap (v) \ra e$ as
$l(\ap ) \ra \infty$, $K\subset U$ for any compact open subgroup in $V$.  Since $V$
is totally disconnected, $K=(e)$.  Thus, $\mu$ has no idempotent factors. 
\eo

\bt\label{thm1}
Let $T_1, T_2 , \cdots , T_m$ be contractions on a finite-dimensional vector space 
$V$ over real or ${\mathbb Q}_p$ 
such that $T_iT_j = T_jT_i$ for all $1\leq i,j\leq m$.  Suppose 
$\mu = T_1(\mu ) *T_2 (\mu )*\cdots * T_m (\mu ) \in {\cal P}(V)$.  Then $\mu$ is
infinitely divisible without idempotent factors.  In particular, autophage measures
on $V$ are infinitely divisible measures without idempotent factors.
\et

\bo
We will treat the case when $m=2$, the general case may be proved in a similar
fashion.  Let $k \geq 1$ be such that $||T^k|| <1$ and $||S^k|| <1$.  
Since $T$ and $S$ commute with each other, we get that 
$\mu = \prod _{i=0} ^{2k} T^{2k-i} S^i (\mu)$.  It is easy to see that
$||T^{2k-i}S^i|| <1$ for any $ 0\leq i \leq 2k$.  By Proposition \ref{pr1}, we get
that $\mu$ is an infinitely divisible measure without idempotent factors.
\eo
\end{section}

\begin{section}{Absolute continuity}

We now introduce a nondegeneracy condition which will be needed for absolute
continuity.  We say that a measure $\mu$ on a finite-dimensional vector space over 
real or ${\mathbb Q}_p$ is {\it full} if $\mu$ is not supported on a coset of a
proper subspace.  We now prove a result on the Fourier transform of full measures. 

\bl\label{lf}
Let $\mu$ be a full measure on a finite-dimensional vector space $V$ over real or
${\mathbb Q}_p$ and $\hat \mu$ be its Fourier transform.  Suppose 
$\mu = T(\mu )* S(\mu )$ for some contractions $T$ and $S$.  Also, if $T$ and $S$
are strict-contractions or $T$ and $S$ commute with each other.  Then
$|\hat \mu (v)|\not =1$ for all $v \in V\setminus (e)$.  
\el

\bo
Here we identify the dual of $V$ with itself.  
Let $\lam = \mu *\check \mu$.  Then it is easy to see that 
\be
\item $\mu $ is full if and only if $\lam$ is not contained in a proper subspace, 

\item for any $v \in V\setminus (e)$, $|\hat \mu (v)| =1$ if and only if $\hat \lam
(v) =1$.
\ee
In view of this it is enough to show that $\hat\lam (v) \not =1$ for any $v \in
\setminus (e)$.  Suppose $\hat\lam (v) =1$ for some $v\not = e$.  Then the closed
subgroup, say $H$, generated by the support of $\lam$ is a proper subgroup (in fact,
$\lam$ is supported on the annihilator of the subgroup of all $v \in V$ such that
$\hat \lam (v)=1 = |\hat \lam (v)|$).  
By assumption, there exists contractions $T$ and $S$ such that $\lam = T(\lam )*
S(\lam )$.  This implies
that $T(H) \subset H$, since $e$ is in the support of $\lam$.  Let $H^0$ be the
maximal vector subspace contained in $H$.  Then $T(H^0)\subset H^0$.  Since $T$ is
an automorphism, $T(H^0) = H^0$.  Similarly, we claim that $S(H^0)=H^0$.  Let
$\tilde \lam$ denotes image of $\lam$ in $V/H^0$ and $\tilde T$ and $\tilde S$ be
the factor automorphisms on $V/H^0$ corresponding to $T$ and $S$.  We now claim that
$\tilde \lam = \delta _e$.  Suppose $V$ is a real vector space.  Then $H/H^0$
is discrete and hence countable.  By Proposition \ref{prp0}, $\lam$ degenerates and
hence $\lam = \delta _e$.  Suppose $V$ is a $p$-adic vector space.  Then $H/H^0 =K$
say, is a compact subgroup.  Since $\tilde \lambda$ is autophage, we have 
$\tilde T(K) \tilde S (K) =K$.  This implies that $\omega _K = T(\omega _K) 
S(\omega _K)$.  By Proposition \ref{pr1} or by Theorem \ref{thm1}, $K=(e)$, that is  
$\lambda =\delta _e$.  This implies that $\mu$ is supported on a coset of a proper
subspace $H^0$.
\eo

\bp\label{prp2}
Let $T_1, T_2, \cdots ,T_m$ be a finite-set of strict-contractions on a
finite-dimensional vector space $V$ over real or ${\mathbb Q}_p$.  Suppose
$\mu$ is a full probability measure on $V$ such that 
$\mu = T_1(\mu ) * T_2(\mu ) *\cdots * T_m(\mu)$.  Then $\hat \mu$
is an integrable function.
\ep

\bo
We treat the case when $m=2$, the general case may be proved in a similar fashion.  
Let $T$ and $S$ be two strict-contractions such that $\mu = T(\mu) *S(\mu)$.  Let 
$T^*$ and $S^*$ be the adjoint of $T$ and $S$ respectively.  
Let $t= ||{T^*}^{-1}||^{-1}$ and 
$s = ||{S^*}^{-1}||^{-1}$.  Let $0<r <\infty$ be such that $t^r +s ^r =1.$  

Let $\hat \mu$ be the Fourier transform of $\mu$.  By Theorem \ref{thm1} and
Theorem 4.2 of \cite {PRV}, $\hat \mu$ does not vanish.  Thus, the function $g$
defined by $$g(v) = -||v||^{-r} \log (| \hat \mu (v)|)$$ for any 
$v \in V \setminus (e)$ is a continuous function.  Since $\mu$ is full, 
by Lemma \ref{lf}, $|\hat \mu (v) |\not = 1$ for any $v \in V\setminus (e)$.  
Thus, $g$ does not vanish.  

Now for any $v \in V\setminus (e)$,
$$\begin{array}{rcl}
g(v) & = & -||v||^{-r} [\log |\hat \mu (T^*(v))| + \log |\hat \mu (S^*(v))|]\\
& = & ({||T^*(v)||\over ||v||})^r g(T^*(v))+({||S^*(v)||\over ||v||})^r g(S^*(v)) \\
&\geq  & t^r g(T^*(v))+s^rg(S^*(v)) \\
\end{array}$$
since for any automorphism $\ap$ of $V$ and $v \in V\setminus (e)$, $||v|| \leq
||\ap ^{-1}|| ||\ap (v)||$.  Thus, for any $v \not = 0 $ in $V$, we have 
$g(v) \geq {\rm min} \{ g (T^*(v)), g(S^*(v)) \}$.  For each $v \in V\setminus
(e)$, we now define for $n \geq 0$, $\ap _n \in sg(T^*, S^*)$ inductively as
follows: $\ap _0 = I$ and for $n \geq 1$, $\ap _{n}$ to be the 
either $T^*\ap _{n-1}$ or $ S^*\ap _{n-1} $ if 
$g(T^*\ap _{n-1} (v) ) \leq g(S^*\ap _{n-1}(v))$  or $g(S^*\ap _{n-1} (v)) \leq
g(T^*\ap _{n-1}(v))$ respectively.  This implies that 
$$\begin{array}{rcl}
g(\ap _{n-1}(v) ) &\geq & t^rg(T^*\ap _{n-1}(v)) +s^r g(S^*\ap _{n-1}(v))\\
& \geq & {\rm min} \{g(T^*\ap _{n-1}(v) ), g(S^*\ap _{n-1} (v) )\}\\
& = & g(\ap _n(v))\\
\end{array}$$ for all $n \geq 1$.  Thus, for each $v \in V\setminus (e)$, there
exists a sequence $(\ap _n)$ (possibly, depending on $v$) in $sg(T^*,S^*)$ such that
$l(\ap _n ) =n$, $\ap _n = T^* \ap _{n-1} ~~ {\rm or} ~~ S^*\ap _{n-1}$ and
$g(v) \geq g(\ap _n(v))$.

Let $c = {\rm min}\{g(v) \mid {\rm min} \{t,s \} \leq ||v|| \leq 1 \}$.  We now
claim that $0< c \leq g(v)$ for any $||v||>1$.  Since $g$ does not vanish,
$c>0$.  Let $v$ be such that $||v||>1$.  Let $\ap _n \in sg(T^*,S^*)$ be such that
$l(\ap _n )=n$, $\ap _n = T^*\ap _{n-1} ~~ {or} ~~ S^*\ap _{n-1}$ 
and $g(\ap _n (v) ) \leq g(v)$.  
Then we get that $\ap _n (v) \ra 0$ as $n \ra \infty$ and $g(v) \geq g(\ap _n(v))$
for all $n $.  Let $k\in \mathbb N$ be such that $||\ap _k (v) ||\leq 1$ but
$||\ap _{k-1}(v)|| >1 $.   Since the roles of $T$ and $S$ are similar, we may assume
that  $\ap _k =T ^*\ap _{k-1}$.  
So, $$1 < ||\ap _{k-1}  (v)|| \leq  ||{T^*}^{-1}|| ||\ap _k(v)||$$ and hence 
$t=||{T^*}^{-1}||^{-1} \leq ||\ap _k(v)||\leq 1 $.  This shows  that 
$$g(v) \geq c$$ for all $v$ with $||v|| >1$.  
It follows from the definition of $g$ that  
$$|\hat \mu (v) | \leq \exp (-c||v|| ^r)$$ for all $v$ with $||v|| >1$ and for some 
constant $c>0$.   Thus, $\hat \mu$ is an integrable function.
\eo

\bt\label{thm2}
Let $T_1, T_2 , \cdots , T_m$ be contractions on a finite-dimensional vector space 
$V$ over real or ${\mathbb Q}_p$ such that 
$T_i T_j = T_jT_i$ for all $1\leq i,j\leq m$.  
Suppose $\mu = T_1 (\mu) * T_2 (\mu ) * \cdots *T_m(\mu ) \in {\cal P}(V)$ is 
full.  Then $\mu$ is absolutely continuous with a bounded continuous density.  In
particular, a full autophage measure on $V$ is absolutely continuous with a bounded
continuous density.
\et

\br
Any question regarding absolute continuity should be addressed only in the class of
full measures because any non-full measure is obviously not absolutely
continuous.  In view of this the assumption fullness is necessary in
Proposition \ref{prp2} and Theorem \ref{thm2}.
\er

\bo
Arguing as in Theorem \ref{thm1}, we may assume that each $T_i$ is
strict-contraction.  Now by Proposition \ref{prp2}, we get that
$\hat \mu$, the Fourier transform of $\mu$ is integrable.  This shows by inverse
Fourier transform that the measure $\mu$ is absolutely continuous with a bounded
continuous density.  
\eo

We now prove the absolute continuity of certain semistable measures.  

\bc\label{cor1}
Let $\mu$ be a full measure on a finite-dimensional vector space over real or 
${\mathbb Q}_p$.  Suppose $\mu$ is $(T,n)$-semistable.  Then $\mu$ is 
absolutely continuous with a bounded continuous density.  
\ec 

\bo
Suppose $\mu$ is $(T,n)$-semistable, that is, $T(\mu ^n) =\mu$.  Suppose $V$ is
a real vector space, then since $\mu$ is full by \cite {Ha}, we get that $T$ is a 
contraction.  By Theorem \ref{thm2}, $\mu$ is absolutely continuous with a bounded
continuous density.  
Suppose $V$ is a vector space over ${\mathbb Q}_p$.  Then by \cite {Sh},
there exist a $x \in V$ and a continuous convolution semigroup $(\mu _t)$ in ${\cal
P}(V)$ such that  $T(\mu _t) = \mu _{nt}$ for all $t >0$ and $\mu _1 = x\mu$.  It is
enough to prove $\mu _1$ is absolutely continuous.  
By \cite {DS}, there exist a compact subgroup $K$ of $V$ and a convolution semigroup
$\tilde \mu _t$ in ${\cal P}(C(T))$ where $C(T) = \{ x\in V \mid T^n (x) \ra e \}$ 
such that 
$\mu _t = \omega _k *\tilde \mu _t$, $T(\tilde \mu _t) = \tilde \mu _{nt}$ and
$C(T)\cap K = (e)$.  Since $\mu$ is full, $V$ is spanned by $C(T)$ and $K$.  Let
$W$ be the vector space spanned by $K$.  Then $K$ is open in $W$ and 
$C(T)\cap W= (e)$.  Thus, $V$ is the direct product of $C(T)$ and $W$.  
Since $\mu _1 = \omega _K*\tilde \mu _1$, it is enough to show that 
$\tilde \mu _1$ is absolutely continuous.  Since $\mu$ is full, $\mu _1$ is full and
hence $\tilde \mu _1$ is full in $C(T)$.  By Theorem \ref{thm2}, 
$\tilde \mu _1$ is absolutely continuous with a bounded continuous density.  
\eo

\end{section}

\begin{acknowledgement}
I wish to thank ICTP, Trieste, Italy for providing financial support and
all facilities during my stay at ICTP.  I also wish to thank NBHM, India for
providing travel grant to visit ICTP.
\end{acknowledgement}

\vskip 0.25in

\noindent Present address: Mathematics Section, ICTP, Strada
Costiera 11, 34014 Trieste. Italy. 

\noindent e-amil: creraja@ictp.trieste.it 

\vskip 0.25in 
\noindent Permanent address:  Stat-Math Unit, Indian Statistical Institute, 8th Mile
Mysore Road, R. V. College Post, Bangalore - 560 059. India. 

\noindent e-amil: creraja@isibang.ac.in


\begin{thebibliography}{10}

\footnotesize 

\bibitem[1]{DS} S. G. Dani and R. Shah, Contraction subgroups and semistable 
measures on $p$-adic Lie groups, Math. Proc. Cambridge Philos. Soc. 110 (1991), 
299-306.  

\bibitem[2]{Ha} W. Hazod, Remarks on(semi-)stable probabilities, Probability
measures on groups, VII (Oberwolfach, 1983), Lecture Notes in Math., 1064, Springer,
Berlin, 1984,  182-203.

\bibitem[3]{Lu} A. Luczak, Operator semistable probability measures on $R^N$,
Colloq. Math. 45 (1981), 287-300.

\bibitem[4]{Par} K. R. Parthasarathy, Probability measures on metric spaces,
Probability and Mathematical Statistics, No.3, Academic Press, Inc., New
York-London 1967.  

\bibitem[5]{PRV} K. R. Parthasarathy, R. Ranga Rao and S. R. S. Varadhan,
Probability distributions on locally compact abelian groups, Illinois J. math. 7
(1963), 337-369.

\bibitem[6]{RS} C. R. E. Raja and Riddhi Shah, Factorization of $T$-decomposable
measures on groups, Monatsh. Math. 133 (2001), 223-239.

\bibitem[7]{RuSz} I. Z. Rusza and G. J. Szekely, Algebraic Probability Theory,
Wiley Series in Probability in Mathematics and Statistics, John Wiley and Sons,
Chichester-New York-Singapore, 1988.

\bibitem[8]{Sh} R. Shah, Semistable measures and limit theorems on real and $p$-adic
groups, Monatsh. Math. 115 (1993), 191-213.  

\bibitem[9]{Si} E. Siebert, Strongly operator-decomposable probability measures on
separable Banach spaces, Math. Nachr. 154 (1991), 315-326.

\bibitem[10]{Sz} G. J. Szekely, Autophage and allelophage probability distributions,
Probability theory and mathematical statistics with applications (Visegr\'ad, 1985),
205-211.





\end{thebibliography}
\end{document}